\documentclass[11pt, fleqn]{amsart}

\usepackage[square,compress,comma, numbers,sort]{natbib}
\usepackage[colorlinks=true, citecolor=red, linkcolor=blue]{hyperref}
\usepackage{amsfonts,mathtools}
\usepackage{amsmath}
\usepackage{mathrsfs}

\allowdisplaybreaks[4]

\usepackage{amssymb}
\usepackage{color}
\usepackage{float}
\usepackage{eufrak}

\DeclarePairedDelimiterXPP\pk[1]{\mathbb{P}}\{ \}{}{ #1}
\DeclarePairedDelimiterXPP\E[1]{\mathbb{E}}\{ \}{}{	#1}

\usepackage{xparse}

\NewDocumentCommand{\ceil}{s O{} m}{%
	\IfBooleanTF{#1} 
	{\left\lceil#3\right\rceil} 
	{#2\lceil#3#2\rceil} 
}
\NewDocumentCommand{\floor}{s O{} m}{%
	\IfBooleanTF{#1} 
	{\left\lfloor#3\right\rfloor}
	{#2\lfloor#3#2\rfloor}
}

\definecolor{c20}{rgb}{0.,0.7,0.}
\definecolor{c30}{rgb}{0.,0.,1.}
\definecolor{c40}{rgb}{1,0.1,0.7}
\definecolor{c50}{rgb}{1,0,0}
\definecolor{c60}{rgb}{1,0.9,0.1}
\definecolor{c70}{rgb}{0.50,1.00,0.00}

\numberwithin{equation}{section}
\newtheorem{theo}{Theorem}[section]
\newtheorem{sat}[theo]{Proposition}
\newtheorem{de}[theo]{Definition}
\newtheorem{lem}[theo]{Lemma}

\newtheorem{example}[theo]{Example}
\newtheorem{korr}[theo]{Corollary}
\newtheorem{remark}[theo]{Remark}

\numberwithin{equation}{section}

\newcommand{\prooftheo}[1]{ \textsc{Proof of Theorem} \ref{#1} }

\newcommand{\prooflem}[1]{\textsc{Proof of Lemma} \ref{#1}}

\newcommand{\QED}{\hfill $\Box$}

\newcommand{\COM}[1]{}

\def\IF{\infty}

\newcommand{\R}{\mathbb{R}}

\newcommand{\BQN}{\begin{eqnarray}}
\newcommand{\EQN}{\end{eqnarray}}
\newcommand{\BQNY}{\begin{eqnarray*}}
	\newcommand{\EQNY}{\end{eqnarray*}}

\newcommand{\limit}[1]{\lim_{#1 \to   \infty}}

\newcommand{\kb}[1]{\boldsymbol{#1}}
\newcommand{\vk}[1]{\kb{#1}}

\def\bqny#1{{\begin{eqnarray*} #1 \end{eqnarray*}}}
\def\bqn#1{{\begin{eqnarray} #1 \end{eqnarray}}}

\newcommand{\BS}{\begin{sat}}
	\newcommand{\ES}{\end{sat}}
\newcommand{\BT}{\begin{theo}}
	\newcommand{\ET}{\end{theo}}
\newcommand{\BK}{\begin{korr}}
	\newcommand{\EK}{\end{korr}}

\newcommand{\BEX}{\begin{example}}
	\newcommand{\EEX}{\end{example}}

\newcommand{\BD}{\begin{de}}
	\newcommand{\ED}{\end{de}}

\newcommand{\BIT}{\begin{itemize}}
	\newcommand{\EIT}{\end{itemize}}
\newcommand{\BDI}{\begin{description}}
	\newcommand{\EDI}{\end{description}}

\newcommand{\BRM}{\begin{remark}}
	\newcommand{\ERM}{\end{remark}}

\newcommand{\BEL}{\begin{lem}}
	\newcommand{\EEL}{\end{lem}}

\newcommand{\nelem}[1]{{Lemma \ref{#1}}}

\begin{document}

\title{Parisian ruin probability for two-dimensional Brownian risk model}

\author{Konrad Krystecki}
\address{Konrad Krystecki, Department of Actuarial Science, 
	University of Lausanne,\\
	UNIL-Dorigny, 1015 Lausanne, Switzerland and
	 Mathematical Institute, University of Wroc\l aw, pl. Grunwaldzki 2/4, 50-384 Wroc\l aw, Poland
}
\email{Konrad.Krystecki@unil.ch}

\bigskip

\date{\today}
 \maketitle

 {\bf Abstract:} Let $(W_1(s), W_2(t)), s,t\ge 0$ be a bivariate Brownian motion with standard Brownian motion marginals and constant correlation $\rho \in (-1,1).$ Parisian ruin is defined as a classical ruin that happens over an extended period of time, the so-called time-in-red. We derive exact asymptotics for the non-simultaneous Parisian ruin of the company conditioned on the event of non-simultaneous ruin happening. We are interested in finding asymptotics of such problem as $u \to \IF$ and with the length of time-in-red being of order $\frac{1}{u^2},$ where $u$ represents initial capital for the companies. Approximation of this problem is of interest for the analysis of Parisian ruin probability in bivariate Brownian risk model, which is a standard way of defining prolonged ruin models in the financial markets.

 {\bf Key Words:} multidimensional Brownian motion; Stationary random fields; Extremes;

 {\bf AMS Classification:} Primary 60G15; secondary 60G70

\section{Introduction}
Consider the following Brownian risk model for two portfolios
$$R_i(t)=u_i+c_i t - W_i(t), i=1,2, $$
where the claims $W_i(t), t \ge 0$ are modeled by two standard dependent Brownian motions, initial capitals $u_i > 0$ and premium rates $c_i.$ The following representation of the dependence between the claims has been proposed in \cite{delsing2018asymptotics} and \cite{DIEKER2005}
\BQNY \label{BB}
(W_1(s),W_2(t))=(B_1(s), \rho B_1(t)+ \sqrt{1- \rho^2} B_2(t)), \quad s,t\ge 0,
\EQNY
where $B_1,B_2$ are two independent standard Brownian motions and $\rho \in [-1,1]$. The ruin probability of a single portfolio in the time horizon $[0,T], T>0$ is given by (see e.g., \cite{MandjesKrzys})
\BQNY \label{single}
{ \pi}_T(c_i,u):=\pk*{ \inf_{t\in [0,T]} R_i(t) < 0}&=& \pk*{\sup_{t\in [0,T]} W_i(t)- c_i t> u}\notag \\
&=& \Phi\left(-\frac{ u}{ \sqrt{T}} -{c_i\sqrt{T}}\right)+
e^{-2c_iu}\Phi\left(- \frac{ u}{\sqrt{T}} +{c_i\sqrt{T}}\right)
\EQNY
for  $i=1,2, u\geq 0$, with $\Phi$ the distribution function of an $N(0,1)$ random variable. Since from self-similarity of Brownian motion we have the following equalities in distribution for $c_1'=\frac{c_1}{\sqrt{T}}, u'=\frac{u}{\sqrt{T}}$
$$B(tT)-c_1 t >u \Leftrightarrow \sqrt{T}B(t) - c_1 t > u \Leftrightarrow  B(t) - c_1' t > u',$$
then without loss of generality one can assume $T=1.$ There are at least two different approaches on how to define the extension of the above to the two-dimensional model. Denote $W_i^*(s)= W_i(s)-c_i s, i=1,2.$ Define the simultaneous ruin probability as
$$ \overline\pi_{A,\rho}(c_1,c_2,u, v)= \pk*{ \exists_{ s \in A}:  W_1^*(s) > u,  W_2^*(s) > v}$$
which has been recently studied in \cite{SIM} for $A=[0,1]$. Similarly, define non-simultaneous ruin probability as
$$ \pi_{A \times B, \rho}(c_1,c_2,u,v)= \pk*{ \exists_{ s \in A, t \in B}:  W_1^*(s) > u,  W_2^*(t) > v}$$
which has been studied for the case $A=B=[0,1]$ in \cite{DHK20}. In this contribution we focus on an extensions of the non-simultaneous results of ruin for two-dimensional risk portfolios. In \cite{loeffen2013} Loeffen, Czarna and Palmowski studied the so-called Parisian ruin of a single portfolio, which is defined as
\BQNY \label{singleParisian}
\mathcal{P}^{*}_{A,H(u)}(c,u) := \pk*{\exists_{s' \in A} \forall_{s \in [s',s'+H(u)]} W^*(s)>u},
\EQNY
for some $H(u) \ge 0$ and $A=[0,T].$ This model defines the concept of the ruin as crossing the barrier over the extended period of time, the so-called time in red. It seems more natural than the classical ruin approach, since it allows for easier practical investigations whether the ruin has occurred. This model has also been studied for various sets $A$ and various processes in many other contributions, e.g. \cite{dkebicki2015parisian}, \cite{dkebicki2016parisian}, \cite{dassios2008parisian}. To analyse the model in two-dimensional framework we use the following definition of the ruin probability
\BQNY \label{doubleParisian}
\mathcal{P}^{*}_{A \times B,\vk H(u)}(c_1,c_2,u,v) := \pk*{\exists_{s' \in A,t' \in B} \forall_{s \in [s',s'+H_1(u)]}\forall_{t \in [t',t'+H_2(u)]} W_1^*(s)>u,W_2^*(t)>v},
\EQNY
for some $H_1(u),H_2(u) \ge 0$ and intervals $A, B.$ We refer to \cite{dassios2011semi}, where one can find an application of Parisian ruin to actuarial risk theory, where $R_i$ is treated as a surplus process of an insurance company with initial capital $u_i.$ For more general intervals $A,B$ we have the following comparison between Parisian and classical ruin
\BQN \label{comparison}
\pi_{A \times B,\rho}(c_1,c_2,u,au) \ge \mathcal{P}^{*}_{A \times B,\vk H(u)}(c_1,c_2,u,au).
\EQN
Since Parisian ruin probability cannot be determined explicitely for general Gaussian risks, our aim is to investigate the asymptotic behaviour of the Parisian ruin conditioned on the classical ruin occurring, for which the results are known. Hence we calculate
\BQNY
\lefteqn{\mathscr{P}^*_{[0,1]^2,\vk H(u)}(c_1,c_2,u,au):=}&&\\
&&\pk*{\exists_{s',t' \in [0,1]} \forall_{s \in [s',s'+H_1(u)]}\forall_{t \in [t',t'+H_2(u)]}\begin{array}{ccc} W_1(s)-c_1s>u \\ W_2(t)-c_2t>au \end{array}\Bigg{|}\exists_{v,w \in [0,1]}\begin{array}{ccc} W_1(v)-c_1v>u \\ W_2(w)-c_2w>au \end{array}}
\EQNY
for $u \to \IF$ and also find for which $H(u)$ we have that for some $C>0$
$$\lim_{u\to \IF}\mathscr{P}^*_{[0,1]^2,\vk H(u)}(c_1,c_2,u,au)=C.$$
We prove that the above is true for $H(u)=(\frac{S_1}{u^2},\frac{S_2}{u^2}):=\frac{(S_1,S_2)}{u^2}$ for some $S_1,S_2>0.$ To simplify notation we denote
$$\mathscr{P}_{S_1,S_2}(c_1,c_2,u,au) := \mathscr{P}^{*}_{[0,1]^2,\frac{(S_1,S_2)}{u^2}}(c_1,c_2,u,au)$$
and similarly
$$\mathcal{P}_{S_1,S_2}(c_1,c_2,u,au) := \mathcal{P}^{*}_{[0,1]^2,\frac{(S_1,S_2)}{u^2}}(c_1,c_2,u,au).$$
Arbitrary choice of $H(u) = \frac{(S_1,S_2)}{u^2}$ is closely connected to the length of the intervals with comparable variance for the Brownian motion (see \cite{pickands1969}). For the choice of $H(u)=o(\frac{1}{u^2})$ following the same line of proof we have that
$$\lim_{u\to \IF}\mathscr{P}^*_{[0,1]^2,\vk H(u)}(c_1,c_2,u,au)=1.$$
On the other hand, if we choose $H(u)$ such that $u^2H(u) \rightarrow \IF, H(u)<1,$ then the methods employed in this contribution are not sufficient and the asymptotics are of different order, even in the one-dimensional setting.

\section{Main results}
Based on the relation between $a$ and $\rho,$ either both of the coordinates impact the asymptotics, or one of the coordinates is negligible (up to a constant). We begin with cases where one of the coordinates dominates the other one and hence the results can be derived from one-dimensional models. Denote by $\Psi$ the survival function of a standard Normal random variable and by $\phi_{t^*}$ the probability density function of $(W_1(1),W_2(t^*)).$
Let $C_{\mathcal{P}}=\mathbb{E}\left\{\exp \left(\sup\limits_{t \ge 0} \inf\limits_{s \in [0,\frac{S_1}{2}]}\sqrt{2}B(t-s)-2|t-s|\mathbf{1}(t>s) \right)\right\},$ which by \cite{dkebicki2016parisian}[Cor 3.5] is positive and finite.
\BT \label{simpleParisian} If $a \le \rho, $ then
\bqny{
\lim_{u \to \IF}\mathscr{P}_{S_1,S_2}(c_1,c_2,u,au) &=& \frac{C_{\mathcal{P}}}{2}.
}
\ET
Our next results are separated into different cases, based on a relative relation between $\rho$ and $A_a= \frac{1}{4a}(1-\sqrt{8a^2+1}).$ Function $A_a$ has been found by analytical calculations. Heuristically, when $\rho<0$ is relatively big compared to $a$ (in terms of absolute value), then it is less likely that the ruin will occur simultaneously and the asymptotics should be significantly different than the ones that have been discovered for simultaneous ruin in \cite{NIKOLAI}. \\
Denote $t_*=\frac{a}{\rho(2a\rho-1)}$ and introduce the following notation for the one-dimensional constants
$$\mathcal{P}(w_1,w_2,f(u)):=\int_{\R}\pk*{ \exists_{ s'\in [0,\infty)} \forall_{s \in [s',s'+f(u)[}: B(s) - w_1 s> x} e^{ w_2 x} dx,$$
$$\mathcal{H}(w_1,w_2,f(u)):=\lim_{\Delta \to \infty}\int_{\R} \frac{1}{\Delta}\pk*{\exists_{ t'\in [0,\Delta]} \forall_{t \in [t',t'+f(u)]} B(t) - w_1t> x} e^{ w_2x} dx,$$
$$\mathcal{R}_{S_1,S_2}=\int_{\R^2}\pk*{\exists_{s',t' \in [0,\infty)} \forall_{s \in [s',s'+S_1], t \in [t',t'+S_2]}:
    \begin{array}{ccc}
	W_1(s)-s>x \\
    W_2(t)-at>y
	\end{array}} e^{\lambda_1 x + \lambda_2 y} dxdy \in (0,\IF).$$
In each particular case, finintess and positivity of $\mathcal{P}$ and $\mathcal{H}$ has been proven in \nelem{1dFiniteParisian}.
\BT \label{MainParisian} Let   $\rho \in (-1,1)$  and $a\in  (\max(0,\rho), 1]$ be given. \\
(i) If $\rho> A_a $, then
\bqn{ \label{M1}
	\lim_{u \to \IF}\mathscr{P}_{S_1,S_2}(c_1,c_2;u,au) = \frac{\mathcal{R}_{S_1,S_2}}{\mathcal{R}_{0,0}}.
}
\newline
(ii) If $\rho= A_a $ and $a <1$, then
\bqn{ \label{M2} \lim_{u \to \IF}\mathscr{P}_{S_1,S_2}(c_1,c_2;u,au) = \frac{(1-a\rho)\mathcal{P}(\frac{1-a\rho}{1-\rho^2},\frac{1-a\rho}{1-\rho^2},S_1)\mathcal{H}(a,2a,S_2)}{2a(1-\rho^2)}.
}
(iii) If $\rho=A_a, a=1$, then
\bqn{ \label{M4}
	\lim_{u \to \IF}\mathscr{P}_{S_1,S_2}(c_1,c_2;u,au) =  \frac{C_{4,1} C'_{4,1}+C_{4,2} C'_{4,2}}{C_4},
}
where $C_{4,1}=\mathcal{P}(2,2,S_1) \mathcal{H}(1,2,S_2),C_{4,2}=\mathcal{P}(2,2,S_2) \mathcal{H}(1,2,S_1)$ and
$$C'_{4,1}=
\begin{cases}
e^{-2\frac{(\frac{1}{2}c_1+c_2)^2}{3}} \Phi\left(c_2+\frac{1}{2}c_1\right), &-\frac{1}{2}c_1<c_2 \\
1, & otherwise,
\end{cases}, \quad C'_{4,2}=
\begin{cases}
e^{-2\frac{(\frac{1}{2}c_2+c_1)^2}{3}}  \Phi\left(c_1+\frac{1}{2}c_2\right), & -\frac{1}{2}c_2<c_1\\
1, & otherwise,
\end{cases} $$
$$ C_4=
\begin{cases}
e^{-2\frac{(\frac{1}{2}c_1+c_2)^2}{3}} \Phi\left(c_2+\frac{1}{2}c_1\right)
+e^{-2\frac{(\frac{1}{2}c_2+c_1)^2}{3}}  \Phi\left(c_1+\frac{1}{2}c_2\right), & c_2>\max(-\frac{1}{2}c_1,-2c_1)\\
e^{-2\frac{(\frac{1}{2}c_1+c_2)^2}{3}} \Phi\left(c_2+\frac{1}{2}c_1\right)
+\frac{1}{2}, &-\frac{1}{2}c_1<c_2 \le -2c_1 \\
\frac{1}{2}
+e^{-2\frac{(\frac{1}{2}c_2+c_1)^2}{3}}  \Phi\left(c_1+\frac{1}{2}c_2\right), &-2c_1<c_2 \le -\frac{1}{2}c_1 \\
1 , &c_2\le \min(-\frac{1}{2}c_1,-2c_1).
\end{cases} $$
(iv) If $a<1, \rho< A_a $, then
\bqn{\label{M5} \lim_{u \to \IF}\mathscr{P}_{S_1,S_2}(c_1,c_2;u,au) = -\frac{\mathcal{P}(\frac{1-a\rho}{1-\rho^2 t_*},\frac{1-a\rho}{1-\rho^2 t_*},S_1) \mathcal{H}(\frac{a}{t_*}, \frac{2a}{t_*},S_2)}{2\rho}.}
(v) If $a=1, \rho< A_a $, then
\bqn{\label{M6} \lim_{u \to \IF}\mathscr{P}_{S_1,S_2}(c_1,c_2;u,au) =  -\frac{C_6}{2\rho},} where
 $t_*=\frac{1}{\rho(2\rho-1)}, C_6=\begin{cases} \mathcal{P}(\frac{1-\rho}{1-\rho^2 t_*},\frac{1-\rho}{1-\rho^2 t_*},S_1) \mathcal{H}(\frac{1}{t_*}, \frac{2}{t_*},S_2)&c_1 \le c_2 \\ \mathcal{P}(\frac{1-\rho}{1-\rho^2 t_*},\frac{1-\rho}{1-\rho^2 t_*},S_2) \mathcal{H}(\frac{1}{t_*}, \frac{2}{t_*},S_1),& c_1 > c_2 \end{cases}.$
\ET

\section{Proofs}
We recall that
\BQNY
\mathscr{P}_{S_1,S_2}(c_1,c_2,u,au) &=& \frac{\mathcal{P}_{S_1,S_2}(c_1,c_2,u,au)}{\pi_{[0,1]^2, \rho}(c_1,c_2,u,au)}.
\EQNY
Therefore in the proofs we can also focus on investigating the asymptotic behaviour of \\ $\mathcal{P}_{S_1,S_2}(c_1,c_2,u,au),$ since the asymptotics for $\pi_{[0,1]^2, \rho}(c_1,c_2,u,au)$ has been calculated in \cite{DHK20}.
\subsection{Proof of Theorem \ref{simpleParisian}}
We divide the proof into two parts - $a < \rho$ and $a=\rho,$ since the methods used are quite different. Further define $S_{1,2}=\max(S_1,S_2),$ which will be commonly used notation in both parts of the proof. \\
\underline{Case (i): $a < \rho.$} First note that
$$ \mathcal{P}_{S_1,S_2}(c_1,c_2,u,au) \le  \pk*{\exists_{s' \in [0,1]} \forall_{s \in [s',s'+\frac{S_1}{u^2}]} W_1^*(s)>u}.$$
On the other hand
\BQNY
\lefteqn{\mathcal{P}_{S_1,S_2}(c_1,c_2,u,au)} && \\ &\ge& \pk*{\exists_{t' \in [0,1]}\begin{array}{ccc} \forall_{t \in [t',t'+\frac{S_1}{u^2}]} B_1^*(t)>u,\rho B_1(t)+ \sqrt{1-\rho^2} B_2(t) - c_2t>au\\\forall_{t \in [t'+\frac{S_1}{u^2},t'+\frac{S_{1,2}}{u^2}]} B_1^*(t)>u-\frac{1}{\sqrt{u}},\rho B_1(t)+ \sqrt{1-\rho^2} B_2(t) - c_2t>au \end{array}}\\
&\ge&\pk*{\exists_{t' \in [0,1]}\begin{array}{ccc} \forall_{t \in [t',t'+\frac{S_1}{u^2}]} B_1^*(t)>u,\rho (u +c_1 t)+ \sqrt{1-\rho^2} B_2(t) - c_2t>au\\ \forall_{t \in [t'+\frac{S_1}{u^2},t'+\frac{S_{1,2}}{u^2}]} B_1^*(t)>u-\frac{1}{\sqrt{u}},\rho (u-\frac{1}{\sqrt{u}} +c_1 t)+ \sqrt{1-\rho^2} B_2(t) - c_2t>au \end{array}}\\
&=&\pk*{\exists_{t' \in [0,1]}\begin{array}{ccc} \forall_{t \in [t',t'+\frac{S_1}{u^2}]} B_1^*(t)>u,\sqrt{1-\rho^2} B_2(t)>(a-\rho)u+(c_2-\rho c_1)t\\ \forall_{t \in [t'+\frac{S_1}{u^2},t'+\frac{S_{1,2}}{u^2}]} B_1^*(t)>u-\frac{1}{\sqrt{u}},\sqrt{1-\rho^2} B_2(t)>(a-\rho)u+\frac{\rho}{u^2}+(c_2-\rho c_1)t \end{array}}\\
&\ge&\pk*{\forall_{s \in [0,1]} B_2(s)>\frac{(a-\rho)u+(c_2-\rho c_1)s}{\sqrt{1-\rho^2}},\exists_{t' \in [0,1]}\begin{array}{ccc} \forall_{t \in [t',t'+\frac{S_1}{u^2}]} B_1(t)-c_1t>u\\ \forall_{t \in [t'+\frac{S_1}{u^2},t'+\frac{S_{1,2}}{u^2}]} B_1^*(t)>u-\frac{1}{\sqrt{u}}\end{array}}\\
&=&\pk*{\exists_{t' \in [0,1]}\begin{array}{ccc} \forall_{t \in [t',t'+\frac{S_1}{u^2}]} B_1^*(t)>u\\ \forall_{t \in [t'+\frac{S_1}{u^2},t'+\frac{S_{1,2}}{u^2}]} B_1^*(t)>u-\frac{1}{\sqrt{u}}\end{array}}\pk*{\forall_{s \in [0,1]} B_2(s)>\frac{(a-\rho)u+(c_2-\rho c_1)s}{\sqrt{1-\rho^2}}}.
\EQNY
Since $a<\rho$, we have that
$$\lim_{u \to \IF}\pk*{\forall_{s \in [0,1]}B_2(s)>\frac{(a-\rho)u+(c_2-\rho c_1)s}{\sqrt{1-\rho^2} }} = 1.$$
Further from independence of increments of Brownian motion we have for $B$ a Brownian motion independent of $B_1,B_2$
\BQNY
\lefteqn{\pk*{\exists_{t' \in [0,1]}\begin{array}{ccc} \forall_{t \in [t',t'+\frac{S_1}{u^2}]} W_1^*(t)>u\\ \forall_{t \in [t'+\frac{S_1}{u^2},t'+\frac{S_{1,2}}{u^2}]} W_1^*(t)>u-\frac{1}{\sqrt{u}}\end{array}}} &&\\ &\ge&\pk*{\exists_{t' \in [0,1]}\forall_{t \in [t',t'+\frac{S_1}{u^2}]} W_1^*(t)>u, \forall_{s \in [0,\frac{\max(S_2-S_1,0)}{u^2}]} B(s)+c_1s<\frac{1}{\sqrt{u}}}\\
&=&\pk*{\exists_{t' \in [0,1]}\forall_{t \in [t',t'+\frac{S_1}{u^2}]} W_1^*(t)>u}\pk*{\forall_{s \in [0,\max(S_2-S_1,0)]} B(s)+\frac{c_1s}{u}<\sqrt{u}}.
\EQNY
Finally we have that
\BQNY
\lim_{u \to \IF}\pk*{\forall_{s \in [0,\max(S_2-S_1,0)]} B(s)+\frac{c_1s}{u}<\sqrt{u}}&=&1
\EQNY
and from \cite{dkebicki2016parisian}[Cor 3.5] we have
\BQN \label{1dpar}
\mathcal{P}_{[0,1],\frac{S_1}{u^2}}(c_1,u) &\sim& C_{\mathcal{P}} \Psi(u+c_1)
\EQN
with $C_{\mathcal{P}}=\mathbb{E}\left\{\sup\limits_{t \ge 0} \inf\limits_{s \in [0,\frac{S_1}{2}]}e^{\sqrt{2}B(t-s)-2|t-s|\mathbf{1}(t>s)}\right\} \in (0,\IF).$
Further we recall that from \cite{DHK20}[Thm 2.1] we have that
$$\pi_{[0,1]^2,\rho}(c_1,c_2;u,au)\sim 2 \Psi(u+c_1).$$
This completes the proof of case (i).
\newline
\underline{Case (ii): $a = \rho.$}
Notice that for $\Delta>0$
\BQNY \mathcal{P}_{S_1,S_2}(c_1,c_2;u, au) &\le& \mathcal{P}^{*}_{[1-\frac{\Delta}{u^2},1]^2,\frac{(S_1,S_2)}{u^2}}(c_1,c_2;u, au) + \pi_{[0,1]^2 \setminus [1-\frac{\Delta}{u^2},1]^2,\rho}(c_1,c_2;u, au).
\EQNY
Denote $\overline{\Delta}(u)=[1-\frac{1}{\sqrt{u}},1].$ Then we have that as $u \to \infty$
\BQNY
\lefteqn{\mathcal{P}^{*}_{[1-\frac{\Delta}{u^2},1]^2,\frac{(S_1,S_2)}{u^2}}(c_1,c_2;u, au)} && \\
 &\le& \pk*{\exists_{s',t' \in \overline{\Delta}(u)} \forall_{(s,t) \in [s',s'+\frac{S_1}{u^2}]\times[t',t'+\frac{S_2}{u^2}]} W_1^*(s)>u,W_2^*(t)>au, \forall_{v \in \overline{\Delta}(u)} W_1^*(v)<u+\frac{1}{\sqrt{u}}}\\
&+& \pk*{\exists_{v \in \overline{\Delta}(u)} W_1^*(v)>u+\frac{1}{\sqrt{u}}}:=\mathbb{P}_1+\mathbb{P}_2
\EQNY
Next observe that
\BQNY
\mathbb{P}_1&\le& \pk*{\exists_{s',t' \in \overline{\Delta}(u)} \forall_{(s,t) \in [s',s'+\frac{S_1}{u^2}]\times[t',t'+\frac{S_2}{u^2}]}\begin{array}{ccc} B_1(s)-c_1 s>u \\ B_2(t)-\frac{c_2-\rho c_1}{\sqrt{1-\rho^2}}t>-\frac{\rho}{\sqrt{u}}\end{array},\forall_{v \in \overline{\Delta}(u)} W_1^*(v)<u+\frac{1}{\sqrt{u}}}\\
&\le& \pk*{\exists_{s' \in \overline{\Delta}(u)} \forall_{s \in [s',s'+\frac{S_1}{u^2}]} B_1(s)-c_1 s>u} \pk*{\exists_{t' \in \overline{\Delta}(u)}\forall_{t \in [t',t'+\frac{S_2}{u^2}]}B_2(t)-\frac{c_2-\rho c_1}{\sqrt{1-\rho^2}}t>0}\\
&\le& \pk*{\exists_{s' \in \overline{\Delta}(u)} \forall_{s \in [s',s'+\frac{S_1}{u^2}]} B_1(s)-c_1 s>u}\pk*{\exists_{t \in \overline{\Delta}(u)}B_2(t)-\frac{c_2-\rho c_1}{\sqrt{1-\rho^2}}t>0}\\
&=& \mathcal{P}^{*}_{[0,1],\frac{S_1}{u^2}}(c_1,u)\Phi\left( \frac{\rho c_1-c_2 }{\sqrt{1-\rho^2}}\right)(1+o(1)), \quad u\to \IF.
\EQNY
Notice that with \cite{dkebicki2016parisian} we have for some $C_1>0$
$$\lim_{u \to \IF}\frac{\mathcal{P}^{*}_{\overline{\Delta}(u),\frac{S_1}{u^2}}(c_1,u)}{\Psi(u+c_1)}=C_1$$
and with \cite{MandjesKrzys} we have for some $C_2>0$
$$\lim_{u \to \IF}\frac{\pi_{\overline{\Delta}(u),\rho}(c_1,u)}{\Psi(u+c_1)}=C_2.$$
Hence
$$\lim_{u \to \IF}\frac{\mathcal{P}^{*}_{\overline{\Delta}(u),\frac{S_1}{u^2}}(c_1,u)}{\pi_{\overline{\Delta}(u),\rho}(c_1,u)}=\frac{C_1}{C_2}.$$
Since from \cite{DHK20} [Thm 2.1] we have
$$\mathbb{P}_2=o\left(\pi_{\overline{\Delta}(u),\rho}(c_1,u)\right), \quad \pi_{[0,1]^2 \setminus \overline{\Delta}(u)^2,\rho}(c_1,c_2;u, au)=o\left(\pi_{\overline{\Delta}(u),\rho}(c_1,u)\right),$$
hence as $u \to \IF$
\BQNY \mathcal{P}_{S_1,S_2}(c_1,c_2;u, au) &\le& \Phi\left( \frac{\rho c_1-c_2 }{\sqrt{1-\rho^2}}\right)\mathcal{P}^{*}_{[0,1],\frac{S_1}{u^2}}(c_1,u).
\EQNY
Finally, following calculations from case (i) we have that for $B$ a Brownian motion independent of $B_1,B_2$
\BQNY
\lefteqn{\mathcal{P}^{*}_{\overline{\Delta}(u)^2,\frac{(S_1,S_2)}{u^2}}(c_1,c_2;u, au)} && \\
 &\ge&\pk*{\exists_{t' \in \overline{\Delta}(u)}\begin{array}{ccc} \forall_{t \in [t',t'+\frac{S_1}{u^2}]} W_1^*(t)>u,\sqrt{1-\rho^2} B_2(t)>(c_2-\rho c_1)t\\ \forall_{t \in [t'+\frac{S_1}{u^2},t'+\frac{S_{1,2}}{u^2}]} W_1^*(t)>u-\frac{1}{\sqrt{u}},\sqrt{1-\rho^2} B_2(t)>\frac{\rho}{u^2}+(c_2-\rho c_1)t \end{array}}\\
&\ge& \pk*{\exists_{t' \in \overline{\Delta}(u)}\begin{array}{ccc} \forall_{t \in [t',t'+\frac{S_1}{u^2}]} W_1^*(t)>u\\ \forall_{t \in [t'+\frac{S_1}{u^2},t'+\frac{S_{1,2}}{u^2}]} W_1^*(t)>u-\frac{1}{\sqrt{u}}\end{array}}\\
&&\times \pk*{\forall_{t \in \overline{\Delta}(u)}\sqrt{1-\rho^2} B_2(t)>(c_2-\rho c_1)t}\\
&\ge& \pk*{\exists_{t' \in [0,1]}\forall_{t \in [t',t'+\frac{S_1}{u^2}]} W_1^*(t)>u}\pk*{\forall_{s \in [0,\max(S_2-S_1,0)]} B(s)+\frac{c_1s}{u}<\sqrt{u}}\\
&&\times \pk*{\exists_{t' \in \overline{\Delta}(u)}\forall_{t \in (t',t'+\frac{S_{1,2}}{u^2})}\sqrt{1-\rho^2} B_2(t)>(c_2-\rho c_1)t}.
\EQNY
Further we have
\BQNY
\pk*{\forall_{t \in \overline{\Delta}(u)}\sqrt{1-\rho^2} B_2(t)>(c_2-\rho c_1)t} &\le& \pk*{\exists_{t \in \overline{\Delta}(u)}\sqrt{1-\rho^2} B_2(t)>(c_2-\rho c_1)t}.
\EQNY
On the other hand with self-similarity and independence of increments of Brownian motion we have that for $B,\hat{B}$ Brownian motions independent of $B_1,B_2$
\BQNY
\lefteqn{\frac{\pk*{\forall_{t \in \overline{\Delta}(u)}\sqrt{1-\rho^2} B_2(t)>(c_2-\rho c_1)t}}{\pk*{\sqrt{1-\rho^2} B_2(1-\frac{1}{\sqrt{u}})>(c_2-\rho c_1)(1-\frac{1}{\sqrt{u}})+\frac{1}{\sqrt[8]{u}}}}} && \\
&\ge& \pk*{\forall_{t \in \overline{\Delta}(u)}\sqrt{1-\rho^2} B_2(t)>(c_2-\rho c_1)t \Bigg{|}\sqrt{1-\rho^2}B_2(1-\frac{1}{\sqrt{u}})>(c_2-\rho c_1)(1-\frac{1}{\sqrt{u}})+\frac{1}{\sqrt[8]{u}}}\\
&\ge&\pk*{\forall_{s \in [0,\frac{1}{\sqrt{u}}]} B(s)-(c_2-\rho c_1)s<\frac{1}{\sqrt[8]{u}} }\\
&=&\pk*{\forall_{s \in [0,1]} \frac{1}{\sqrt[4]{u}}B(s)-\frac{1}{\sqrt{u}}(c_2-\rho c_1)s<\frac{1}{\sqrt[8]{u}}}\\
&=&\pk*{\forall_{s \in [0,1]} B(s)-\frac{1}{\sqrt[4]{u}}(c_2-\rho c_1)s<\sqrt[4]{u}} \sim 1
\EQNY
Finally
\BQNY
\lefteqn{\lim_{u \to \IF}\pk*{\exists_{t \in \overline{\Delta}(u)}\sqrt{1-\rho^2} B_2(t)>(c_2-\rho c_1)t}}&&\\
&=&\lim_{u \to \IF}\pk*{\sqrt{1-\rho^2} B_2(1-\frac{1}{\sqrt{u}})>(c_2-\rho c_1)(1-\frac{1}{\sqrt{u}})+\frac{1}{\sqrt[8]{u}}}= \Phi\left( \frac{\rho c_1-c_2 }{\sqrt{1-\rho^2}}\right).
\EQNY
Hence the claim follows from \eqref{1dpar} and from \cite{DHK20}[Thm 2.1], which gives
$$\pi_{[0,1]^2,\rho}(c_1,c_2;u,au)\sim 2 \Phi\left( \frac{\rho c_1-c_2 }{\sqrt{1-\rho^2}}\right) \Psi(u+c_1).$$
\QED
\subsection{Proof of Theorem \ref{MainParisian}}
We again recall that
\BQNY
\mathscr{P}_{S_1,S_2}(c_1,c_2,u,au) &=& \frac{\mathcal{P}_{S_1,S_2}(c_1,c_2,u,au)}{\pi_{[0,1]^2, \rho}(c_1,c_2,u,au)}
\EQNY
and as in the proof of the Theorem \ref{simpleParisian} we focus on investigating the asymptotic behaviour of \\ $\mathcal{P}_{S_1,S_2}(c_1,c_2,u,au).$ Before we begin the proof we need few technical lemmas. First let
$$\Sigma_{ s,t}= \begin{pmatrix}
s&  \rho \min(s,t) \\
\rho \min(s,t) &  t
\end{pmatrix} $$
be the covariance matrix of $(W_1(s),W_2(t))$. In \cite{DHK20} it was noted that the drift has a significant impact on the optimization problem that was used to determine asymptotics for the classical ruin. We denote below for $\vk a = (1+\frac{c_1 s}{u},a+\frac{c_2t}{u})^\top$
 $$ q_{\vk a}(s,t):= \vk a ^\top \Sigma^{-1}_{s,t} \vk a ,\quad \vk{b}(s,t):=  \Sigma^{-1}_{s,t} \vk a$$
 and set
\bqn{
	 q_{\vk a}^*(s,t)= \min_{ \vk x \ge \vk a} q_{\vk x}(s,t), \quad  q_{\vk a}^*= \min_{s,t \in [0,1]}  q_{\vk a}^*(s,t).
	}
Note that for $a > \rho$ and large enough $u$ we have $\vk b(s,t) \sim (\frac{t-a\rho\min(s,t)}{st-\rho^2(\min(s,t))^2},\frac{as-\rho\min(s,t)}{st-\rho^2(\min(s,t))^2}) > \vk 0.$ From \cite{DEBICKIKOSINSKI} we have that for any $s,t$ positive the following logarythmic asymptotics occurs
\bqn{ \label{EXA2}
 \limit{u} \frac{1}{u^2} 	\log \pk{\exists_{s,t \in [0,1]} W_1^*(s)> u, W_2^*(t) > au} =  - \frac{q_{\vk a}^*(s,t)}{2} .
}	
Hence we will use the function $q_{\vk a}^*(s,t)$ to reflect the asymptotics of $\pk{ W_1^*(s)> u, W_2^*(t) > au}$. Below we present the main lemma that solves the optimization problem stated above and was first derived in \cite{DHK20}.
\BEL \label{optDrift}
For all large $u$ we have:\\
(i) If $a=1, \rho<-\frac{1}{2}$,  then   $q_{\vk a_u(s,t)}^*(s,t)$ attains its unique local minima on $[0,1]^2$ at
$$(s_u,t_u):= \left(1,\frac{1}{\rho(2\rho-1)+\frac{c_2-\rho c_1}{u}}\right), \quad (\bar{s}_u,\bar{t}_u):= \left(\frac{1}{\rho(2\rho-1)+\frac{c_1-\rho c_2}{u}}, 1\right).$$
(ii) If $a=1, \rho=-\frac{1}{2},$ then
$q_{\vk a_u(s,t)}^*(s,t)$ attains its unique local minima on $[0,1]^2$ at
$$(s_u,t_u):= \left(1,\min(\frac{1}{1+\frac{c_2+2 c_1}{u}},1)\right), \quad (\bar{s}_u,\bar{t}_u):= \left(\min(\frac{1}{1+\frac{c_1+2 c_2}{u}},1), 1\right).$$
(iii) For any other $a \in (\max(0,\rho),1], \rho \in (-1,1),$ $q_{\vk a_u(s,t)}^*(s,t)$ attains its unique minimum on $[0,1]^2$ at
 $$(s_u,t_u):=\begin{cases}
 (1,\frac{a}{\rho(2a\rho-1)+\frac{c_2-\rho c_1}{u}}) , & { if} \  \frac{a}{\rho(2a\rho-1)+\frac{c_2-\rho c_1}{u}} \in [0, 1] \\
 (1,1), & otherwise. \end{cases}
  $$
\EEL

In the rest of the paper we denote
$$t^*:=\lim_{u \to \IF}t_u,$$
where $t_u$ is defined as in \nelem{optDrift}. We moreover recall the comparison between the behaviour of variance in $t_u$ and $t^*$ from \cite{DHK20}[Rem. 3.3].
\BRM \label{32} For $a \in (\max(0,\rho),1]$ and $t_u<1$ as in \nelem{optDrift} we have
$$\varphi_{t_u}(u+c_1 ,au+c_2t_u) \sim e^{-a\frac{(c_1\rho-c_2)^2}{2\rho(1-a\rho)}}\varphi_{t^*}(u+c_1 ,au+c_2t^*), \ u \to \IF.$$
\ERM
Denote $k_{u}=  1- \frac{(k-1)\Delta}{u^2},l_{u}=  t_u- \frac{(l-1)\Delta}{u^2}, u>0, \Delta>0$ and set
$$E_{u,k}=[(k+1)_u, k_u],E_{u,k,l}=E_{u,k} \times E_{u,l}, \quad E= [-\Delta,0] \times [-\Delta,0].$$ Define also $\eta_{u,k,l}(s,t):=(\eta_{1,u,k}(s),\eta_{2,u,l}(t)):=u(W_1(\frac{s}{u^2}+k_u) - W_1(k_u) - c_1\frac{s}{u^2}, W_2(\frac{t}{u^2}+l_u )- W_2(l_u) - c_2\frac{t}{u^2}).$
The following lemma is used to calculate the ruin probability on an interval of size of order $O(\frac{1}{u^2}).$
\BEL \label{PickandsParisian}
Let $\rho \in (-1,1), a \in (\max(0,\rho),1], l,k=O(\frac{u\log(u)}{\Delta}) $ and  $\Delta, S_1,S_2>0 $ be given constants. Then, as $u\to \IF$
\BQNY
\mathcal{P}^{*}_{E_{u,k,l},\frac{(S_1,S_2)}{u^2}}(c_1,c_2,u,au)&\sim &u^{-2} \varphi_{t^*}(u+c_1 ,au+c_2 t_u) I_1(\Delta) e^{-\frac{1}{2}u^2 (q_a^*(k_u,l_u)-q_a^*(1,t_u))},
\EQNY
$I_1(\Delta)= 	\begin{cases}
\int_{\R^2}\pk*{\exists_{s',t' \in [0,\Delta]} \forall_{s \in [s',s'+S_1], t \in [t',t'+S_2]}:
    \begin{array}{ccc}
	W_1(s)-s>x \\
    W_2(t)-at>y
	\end{array}}  e^{\lambda_1 x + \lambda_2 y} dxdy &l_u = k_u\\
\begin{aligned}  \int_{\R^2}\pk*{\exists_{t' \in [0,\Delta]}\forall_{t \in [t',t'+S_2]}:  W_2(t)-\frac{a-\rho}{t^*-\rho^2}t>y} \\
\times \pk*{\exists_{s' \in [0,\Delta]}\forall_{s \in [s',s'+S_1]}: W_1(s)-s>x}  e^{\lambda_1 x + \lambda_2 y} dxdy \end{aligned} &l_u > k_u \\
\begin{aligned}\int_{\R^2}\pk*{\exists_{s' \in [0,\Delta]}\forall_{s \in [s',s'+S_1]}: W_1(s)-\frac{1-a\rho}{1-\rho^2t^*}s>x} \\
\times \pk*{\exists_{t' \in [0,\Delta]}\forall_{t \in [t',t'+S_2]}:  W_2(t)-\frac{a}{t^*}t>y}  e^{\lambda_1 x + \lambda_2 y} \end{aligned} &l_u < k_u
\end{cases},$
\\
and $\lambda_1 = \begin{cases}
\frac{1}{t^*}\frac{1-a\rho}{1-\rho^2} &l_u = k_u\\
\frac{t^*-a\rho}{t^*-\rho^2}, &l_u > k_u \\
\frac{1-a\rho}{1-\rho^2t^*}, &l_u < k_u \\
\end{cases}$, $\lambda_2 = \begin{cases}
\frac{1}{t^*}\frac{a-\rho}{1-\rho^2}, &l_u = k_u\\
\frac{a-\rho}{t^*-\rho^2 }, &l_u > k_u \\
\frac{a-\rho t^*}{t^* -\rho^2(t^*)^2}, &l_u < k_u \\
\end{cases}.$
\newline
Additionally
\scriptsize{
\begin{equation}
\label{Int1}
 \lim_{u \to \infty}\sup_{l,k = O(u\log u)}\int_{\R^2}\pk*{\exists_{(s',t') \in E} \forall_{s \in [s',s'+S_1], t \in [t',t'+S_2]}:
    \begin{array}{ccc}
	\eta_{u,k,l}(s,t)>(x,y)
	\end{array}
    \Bigg{|}
    \begin{array}{ccc}
	W_1^*(k_u )= u - \frac{x}{u}  \\
	W_2^*(l_u )= au - \frac{y}{u}
	\end{array} }  e^{\lambda_1 x + \lambda_2 y} dxdy\\
< \infty.
\end{equation}}
\EEL

\BRM
Limit \eqref{Int1} is necessary so that dominated converence theorem can be used while proving \nelem{PickandsParisian}.
\ERM

The following lemmas are used to show that the constants $I_1$ in \nelem{PickandsParisian} are positive and finite.
\BEL \label{1dFiniteParisian}
i) For any $b, c>0, S \ge 0$ such that  $2b>c$  we have
$$
\int_{\R} \pk{ \exists_{ t' \ge 0} \forall_{t \in [t',t'+S]} W(t) - bt> x} e^{ cx} dx \in (0,\infty).
$$
ii) For any $b>0, S \ge 0$
$$
\lim_{T\to\infty}\frac{1}{T}\int_{\R} \pk{ \exists_{ t' \in [0,T]} \forall_{t \in [t',t'+S]} W(t) - bt> x} e^{ 2bx} dx \in (0,\infty).
$$
\EEL
\prooflem{1dFiniteParisian} \\
Ad i) We have that
\BQNY
\int_{\R} \pk{ \exists_{ t' \ge 0} \forall_{t \in [t',t'+S]} W(t) - bt> x} e^{ cx} dx &=& \int_{\R} \pk{ \exists_{ t' \ge 0} \forall_{t \in [0,S]} W(t+t') - b(t+t')> x} e^{ cx}dx\\
&=& \frac{1}{b}\int_{\R} \pk{ \exists_{ t' \ge 0} \forall_{t \in [0,S]} bW(t+t') - b^2(t+t')> x} e^{ \frac{c}{b}x}dx\\
&=& \frac{1}{b}\int_{\R} \pk{ \exists_{ t' \ge 0} \forall_{t \in [0,\frac{S}{b^2}]} W(t+t') - (t+t')> x} e^{ \frac{c}{b}x}dx
\EQNY
This constant is the same as in \cite{dkebicki2016parisian}(8) for $T=\frac{S}{b^2}, \alpha=1, \beta=1, b_1=b_2=\frac{2b-c}{2b}.$ \\
Ad ii) Notice that
\BQNY
\lefteqn{\lim_{T\to\infty}\frac{1}{T}\int_{\R} \pk{ \exists_{ t' \in [0,T]} \forall_{t \in [t',t'+S]} W(t) - bt> x} e^{ 2bx} dx} &&\\
 &=& \lim_{T\to\infty}\frac{1}{T}\int_{\R} \pk{ \exists_{ t' \in [0,T]} \forall_{t \in [0,S]} W(t+t') - b(t+t')> x} e^{ 2bx} dx \\
 &=& \lim_{T\to\infty}\frac{1}{bT}\int_{\R} \pk{ \exists_{ t' \in [0,T]} \forall_{t \in [0,S]} bW(t+t') - b^2(t+t')> x} e^{ 2x} dx \\
 &=& \lim_{T\to\infty}\frac{1}{bT}\int_{\R} \pk{ \exists_{ t' \in [0,\frac{T}{b^2}]} \forall_{t \in [0,\frac{S}{b^2}]} W(t+t') - (t+t')> x} e^{ 2x} dx.
\EQNY
This constant is the same as in \cite{dkebicki2015parisian}(2.5) for $\alpha=1$ and $S=\frac{T}{b^2},T=\frac{S}{b^2}.$
\QED

\BEL \label {2dFiniteParisian} Take any $a > \max(0,\rho), S_1,S_2 \ge 0.$ Then
$$
\int_{\R^2}\pk*{\exists_{s',t' \in \R_+} \forall_{s \in [s',s'+S_1], t \in [t',t'+S_2]}:
    \begin{array}{ccc}
	W_1(s)-s>x \\
    W_2(t)-at>y
	\end{array}}  e^{\lambda_1 x + \lambda_2 y} dxdy \in (0,\infty),
$$
where $\lambda_1 = \frac{1-a\rho}{1-\rho^2}, \lambda_2 = \frac{a-\rho}{1-\rho^2}.$
\EEL
\prooflem{2dFiniteParisian}
Positivity comes of the constant from the fact that the function that we integrate is positive on a set of positive mass. Finitness of the constant follows straightforwardly from
\BQNY
\lefteqn{\int_{\R^2}\pk*{\exists_{s',t' \in \R_+} \forall_{s \in [s',s'+S_1], t \in [t',t'+S_2]}:
    \begin{array}{ccc}
	W_1(s)-s>x \\
    W_2(t)-at>y
	\end{array}}  e^{\lambda_1 x + \lambda_2 y} dxdy } && \\
&\le&\int_{\R^2}\pk*{\exists_{s',t' \in \R_+}:
    \begin{array}{ccc}
	W_1(s')-s'>x \\
    W_2(t')-at'>y
	\end{array}}  e^{\lambda_1 x + \lambda_2 y} dxdy
\EQNY
and \cite{DHK20} [Lemma 3.6].
\QED

\prooftheo{MainParisian}
We split the proof into several cases which depend on the behaviour of the variance and the optimization point we get from \nelem{optDrift}. Let next
$$N_u:=\floor{\frac{u\log(u)}{\Delta}}, \quad K_u^{(1)}=\frac{(c_2-c_1\rho)u}{\Delta}, \quad K_u^{(2)}=\frac{(c_1-c_2\rho)u}{\Delta},$$
$$E_{u,m}^1:=[(m+1)_u, m_u],\quad E_{u,j}^2:=[(j+1)_u, j_u],$$
where $m_{u}=  1- \frac{(m-1)\Delta}{u^2}, j_{u}=  t_u- \frac{(j-1)\Delta}{u^2}$. For $\Delta>0$ we have
$$\mathcal{P}_{S_1,S_2}(c_1,c_2;u, au) \ge \mathcal{P}^{*}_{F_u,\frac{(S_1,S_2)}{u^2}}(c_1,c_2;u, au).$$
On the other hand
$$
\mathcal{P}_{S_1,S_2}(c_1,c_2;u, au) \le \mathcal{P}^{*}_{F_u,\frac{(S_1,S_2)}{u^2}}(c_1,c_2;u, au)+\pi_{[0,1]^2 \setminus F_u}(c_1,c_2;u, au).
$$
We have
\BQNY
\mathscr{P}_{S_1,S_2}(c_1,c_2;u, au) &\le& \frac{\mathcal{P}^{*}_{F_u,\frac{(S_1,S_2)}{u^2}}(c_1,c_2;u, au)+\pi_{[0,1]^2 \setminus F_u}(c_1,c_2;u, au)}{\pi_{[0,1]^2}(c_1,c_2;u, au)}.
\EQNY
Since from \cite{DHK20}[Thm 2.2] we have that
$$\lim_{u \to \IF}\frac{\pi_{[0,1]^2 \setminus F_u}(c_1,c_2;u, au)}{\pi_{[0,1]^2}(c_1,c_2;u, au)}=0$$
therefore
$$\mathscr{P}_{S_1,S_2}(c_1,c_2;u, au) \sim \frac{\mathcal{P}^{*}_{F_u,\frac{(S_1,S_2)}{u^2}}(c_1,c_2;u, au)}{\pi_{[0,1]^2}(c_1,c_2;u, au)}, $$
where $F_u$ is case dependant.
\newline
\underline{Case (i): $\rho>\frac{1}{4a}(1-\sqrt{8a^2+1}).$} According to \nelem{optDrift} $t^*=t_u=1$. From \cite{DHK20}[Thm 2.2, case (i)] we have $F_u:=E_{u,1}^2.$
Using \nelem{PickandsParisian} and \nelem{2dFiniteParisian} and taking $u \to \infty$ and then $\Delta \to \infty$, we get that
\BQNY
\mathcal{P}^{*}_{E_{u,1}^2,\frac{(S_1,S_2)}{u^2}}(c_1,c_2;u, au) &\sim &  I u^{-2} \varphi_{1}(u+c_1 ,au+c_2),
\EQNY
where $I=\int_{\R^2}\pk*{\exists_{s',t' \in [0,\infty)} \forall_{s \in [s',s'+S_1], t \in [t',t'+S_2]}:
    \begin{array}{ccc}
	W_1(s)-s>x \\
    W_2(t)-at>y
	\end{array}} e^{\lambda_1 x + \lambda_2 y} dxdy.$
With that, the proof of case (i) is complete.
\newline
\underline{Case (ii): $\rho=\frac{1}{4a}(1-\sqrt{8a^2+1}).$} We split this case into two subcases since the behaviour of the optimizing point is slightly different. First let , $c_2-\rho c_1 \le 0.$ According to \nelem{optDrift} $t^*=t_u=1$. From \cite{DHK20}[Thm 2.2, case (ii)] we have $F_u:=[1-\frac{\Delta}{u^2},1]\times[1-\frac{\log(u)}{u},1-\frac{\Delta}{u^2}].$
Using Bonferroni inequality we have that
\BQN \label{lBP}
\lefteqn{\mathcal{P}^{*}_{F_u,\frac{(S_1,S_2)}{u^2}}(c_1,c_2;u, au)} && \nonumber\\ &\ge& \sum_{l=2}^{N_u}\pk*{\exists_{s' \in E_{u,1}^1, t' \in E_{u,l}^2}\forall_{s \in [s',s'+\frac{S_1}{u^2}]}, t \in [t',t'+\frac{S_2}{u^2}]: W_1^*(s)>u, W_2^*(t)> au  } \nonumber \\
&&- \sum_{l=2}^{N_u}\sum_{m=l+1}^{N_u}\pk*{\exists_{s \in E_{u,1}^1, t_1 \in E_{u,l}^2, t_2 \in E_{u,m}^2}: W_1^*(s)>u, W_2^*(t_1)> au, W_2^*(t_2)> au  } \nonumber\\
&:=&P_{u,\Delta}-D_{u,\Delta}.
\EQN
Further we have
\BQN \label{uBP}
\mathcal{P}^{*}_{F_u,\frac{(S_1,S_2)}{u^2}}(c_1,c_2;u, au) &\le& P_{u,\Delta}+D_{u,\Delta}.
\EQN
From \nelem{PickandsParisian} we have as $u \to \IF$
\BQNY
P_{u,\Delta}& \sim & C_{2,\mathcal{P}}^{(1)}(\Delta) C_{2,\mathcal{P}}^{(2)}(\Delta) u^{-2}\varphi_{t^*}(u+c_1 ,au+c_2) \sum_{l=2}^{N_u} e^{-\frac{1}{2}u^2(q_{\vk a}(k_u,l_u)-q_{\vk a}(1,1))},\\
\EQNY
where
$$C_{2,\mathcal{P}}^{(1)}(\Delta)=\int_{\R}\pk*{ \exists_{ s'\in [0,\Delta]} \forall_{s \in [s',s'+S_1]}: W_1(s) - \frac{1-a\rho}{1-\rho^2}s> x} e^{ \frac{1-a\rho}{1-\rho^2}x} dx$$
and
$$C_{2,\mathcal{P}}^{(2)}(\Delta)=\int_{\R} \pk*{ \exists_{ t'\in [0,\Delta]} \forall_{t \in [t',t'+S_2]}: W_2(t) - at> x} e^{ 2ax} dx.$$
Using Taylor expansion we have that for $k<l$
$$u^2(q_a(k_u,l_u)-q_a(1,1))=\tau_1 (k-1)\Delta+\tau_4 \frac{(l-1)^2\Delta^2}{u^2}+o(\frac{k^2}{u^2})+o(\frac{l^3}{u^4}),$$
where $\tau_1=\frac{(1 - a \rho)^2}{(1 - \rho^2)^2}>0$ and $\tau_4=\frac{\rho^2  - 2 a \rho^3 +  a^2 \rho^2 }{(1 - \rho^2)^2}>0.$
Therefore with \cite{DHK20}[Lem 3.6] we have as $u \to \IF$
\BQNY
P_{u,\Delta} & \sim & C_{2,\mathcal{P}}^{(1)}(\Delta) C_{2,\mathcal{P}}^{(2)}(\Delta) u^{-2}\varphi_{t^*}(u+c_1 ,au+c_2) \sum_{l=2}^{N_u} e^{-\frac{\tau_4}{2}\frac{(l-1)^2\Delta^2}{u^2}}\\
&=& \frac{1}{\sqrt{\tau_4}}C_{2,\mathcal{P}}^{(1)}(\Delta) \frac{C_{2,\mathcal{P}}^{(2)}(\Delta)}{\Delta} u^{-1}\varphi_{t^*}(u+c_1 ,au+c_2) \sum_{l=2}^{N_u} \frac{\sqrt{\tau_4}\Delta}{u} e^{-\frac{\tau_4}{2}\frac{(l-1)^2\Delta^2}{u^2}}\\
&\sim& C_{2,\mathcal{P}}^{(1)}(\Delta) \frac{C_{2,\mathcal{P}}^{(2)}(\Delta)}{\Delta} \frac{\sqrt{\pi}}{\sqrt{2\tau_4}} u^{-1}\varphi_{t^*}(u+c_1 ,au+c_2).
\EQNY
From \nelem{1dFiniteParisian} we have that $\lim_{\Delta \to \IF}C_{2,\mathcal{P}}^{(1)}(\Delta)=C_{2,\mathcal{P}}^{(1)}$ and $\lim_{\Delta \to \IF}\frac{C_{2,\mathcal{P}}^{(2)}(\Delta)}{\Delta}=C_{2,\mathcal{P}}^{(2)}.$ Hence
$$\lim_{\Delta \to \IF}\lim_{u \to \IF}\frac{P_{u,\Delta}}{C_{2,\mathcal{P}}^{(1)}C_{2,\mathcal{P}}^{(2)}\frac{\sqrt{\pi}}{\sqrt{2\tau_4}} u^{-1}\varphi_{t^*}(u+c_1 ,au+c_2)}=1.$$
From the proof of \cite{DHK20}[Theorem 2.2, case (ii)] we have that for $C=C_{2,\mathcal{P}}^{(1)}C_{2,\mathcal{P}}^{(2)}\frac{\sqrt{\pi}}{\sqrt{2\tau_4}}$
\BQN \label{negl1}
\lim_{\Delta \to \IF}\lim_{u \to \IF}\frac{D_{u,\Delta}}{P_{u,\Delta}}&=&\lim_{\Delta \to \IF}\lim_{u \to \IF}\frac{D_{u,\Delta}}{Cu^{-1}\varphi_{t^*}(u+c_1 ,au+c_2)}=0.
\EQN
i.e. the double sum is negligible compared to the single sum.
\newline
Now we consider the case of $c_2-\rho c_1 > 0.$ According to \nelem{optDrift}
there is exactly one minimizer of
$q_{\vk a_u(s,t)}^*(s,t)$  on $[0,1]^2$:  $(s_u,t_u)=\left(1,\frac{a}{\rho(2a\rho-1)+\frac{c_2-\rho c_1}{u}}\right),$ where from $c_1 - \rho c_2 >0$ we obtain that
$$\frac{a}{\rho(2a\rho-1)+\frac{c_2-\rho c_1}{u}} \nearrow 1 $$
as $u \to \IF.$ From \cite{DHK20}[Thm 2.2, case (iii)] we have $F_u:=[1-\frac{\Delta}{u^2},1]\times[t_u-\frac{\log(u)}{u},1-\frac{\Delta}{u^2}].$ The proof follows the path of the proof of case (ii).
Again from \nelem{PickandsParisian} we have
\BQNY
P_{u,\Delta}& \sim & C_{3,\mathcal{P}}^{(1)}(\Delta) C_{3,\mathcal{P}}^{(2)}(\Delta) u^{-2}\varphi_{t^*}(u+c_1 ,au+c_2) \sum_{l=-K_u^{(1)}}^{N_u} e^{-\frac{1}{2}u^2(q_{\vk a}(k_u,l_u)-q_{\vk a}(1,1))},\\
\EQNY
where
$$C_{3,\mathcal{P}}^{(1)}(\Delta)=\int_{\R}\pk*{ \exists_{ s'\in [0,\Delta]} \forall_{s \in [s',s'+S_1]}: W_1(s) - \frac{1-a\rho}{1-\rho^2}s> x} e^{ \frac{1-a\rho}{1-\rho^2}x} dx$$
and
$$C_{3,\mathcal{P}}^{(2)}(\Delta)=\int_{\R} \pk*{ \exists_{ t'\in [0,\Delta]} \forall_{t \in [t',t'+S_2]}: W_2(t) - at> x} e^{ 2ax} dx.$$
Observe that the Taylor expansions mentioned in case (ii) are independent on $c_1,c_2$ and hence remain the same here. Therefore with \cite{DHK20}[Lem 3.6] we have
\BQNY
P_{u,\Delta} & \sim & C_{3,\mathcal{P}}^{(1)}(\Delta) C_{3,\mathcal{P}}^{(2)}(\Delta) u^{-2}\varphi_{t^*}(u+c_1 ,au+c_2) \sum_{l=-K_u^{(1)}}^{N_u} e^{-\frac{\tau_4}{2}\frac{(l-1)^2\Delta^2}{u^2}}\\
&\sim& C_{3,\mathcal{P}}^{(1)}(\Delta) \frac{C_{3,\mathcal{P}}^{(2)}(\Delta)}{\Delta} \frac{\sqrt{2\pi}}{\sqrt{\tau_4}} e^{-a\frac{(c_1\rho-c_2)^2}{2\rho(1-a\rho)}}  \Phi\left(c_2-\rho c_1\right) u^{-1}\varphi_{t^*}(u+c_1 ,au+c_2).
\EQNY
From \nelem{1dFiniteParisian} we have that $\lim_{\Delta \to \IF}C_{3,\mathcal{P}}^{(1)}(\Delta)=C_{3,\mathcal{P}}^{(1)}$ and $\lim_{\Delta \to \IF}\frac{C_{3,\mathcal{P}}^{(2)}(\Delta)}{\Delta}=C_{3,\mathcal{P}}^{(2)}.$ Hence we have that
$$\lim_{\Delta \to \IF}\lim_{u \to \IF}\frac{P_{u,\Delta}}{C_{2,\mathcal{P}}^{(1)}C_{2,\mathcal{P}}^{(2)}\frac{\sqrt{2\pi}}{\sqrt{\tau_4}} e^{-a\frac{(c_1\rho-c_2)^2}{2\rho(1-a\rho)}}  \Phi\left(c_2-\rho c_1\right)  u^{-1}\varphi_{t^*}(u+c_1 ,au+c_2)}=1.$$
Since $D_{u,\Delta}$ is the same as \cite{DHK20}(3.16), we have from \cite{DHK20} that
\BQNY
\lim_{\Delta \to \IF}\lim_{u \to \IF}\frac{D_{u,\Delta}}{Cu^{-1}\varphi_{t^*}(u+c_1 ,au+c_2)}&=&0.
\EQNY
With that the proof of case (ii) is complete.
\newline
\underline{Case (iii): $\rho=-\frac{1}{2}, a=1.$} According to \nelem{optDrift} $t^*=1$.
The proof is analogous to case (ii). We use \eqref{uBP} and \eqref{lBP} with
$$F_u:=[1-\frac{\Delta}{u^2},1]\times[1-\frac{\log(u)}{u},1-\frac{\Delta}{u^2}]\cup [1-\frac{\log(u)}{u},1-\frac{\Delta}{u^2}]\times[1-\frac{\Delta}{u^2},1],$$
\BQNY
P_{u,\Delta}&=&\sum_{l=2}^{N_u}\pk*{\exists_{s' \in E_{u,1}^1, t' \in E_{u,l}^2}\forall_{s \in (s',s'+H)}, t \in (t',t'+H): W_1^*(s)>u, W_2^*(t)> u  }\\
&&+\sum_{k=2}^{N_u}\pk*{\exists_{s' \in E_{u,k}^1, t' \in E_{u,1}^2}\forall_{s \in (s',s'+H)}, t \in (t',t'+H): W_1^*(s)>u, W_2^*(t)> u  }\\
&:=&P_{u,\Delta}^{(1)}+P_{u,\Delta}^{(2)},
\EQNY
\BQNY
D_{u,\Delta}&=&\sum_{l=2}^{N_u}\sum_{m=l+1}^{N_u}\pk*{\exists_{s \in E_{u,1}^1, t_1 \in E_{u,l}^2, t_2 \in E_{u,m}^2}: W_1^*(s)>u, W_2^*(t_1)> au, W_2^*(t_2)> u  }\\
&&+\sum_{k=2}^{N_u}\sum_{m=k+1}^{N_u}\pk*{\exists_{s_1 \in E_{u,k}^1, s_2 \in E_{u,m}^1, t \in E_{u,1}^2}: W_1^*(s_1)>u, W_1^*(s_2)> u, W_2^*(t)> u  }\\
&&+\sum_{k=2}^{N_u}\sum_{l=2}^{N_u}\pk*{\exists_{s_1 \in E_{u,l}^1, s_2 \in E_{u,1}^1, t_1 \in E_{u,1}^2, t_2 \in E_{u,k}^2}:
\begin{array}{ccc}
	W_1^*(s_1)>u \\
    W_1^*(s_2)>u \\
	W_2^*(t_1) > u \\
    W_2^*(t_2) > u
\end{array}
}.
\EQNY
Notice that using calculations from case (ii) and (iii) for $a=1,\rho=-\frac{1}{2}$ we have that
$$P_{u,\Delta}^{(1)}=\begin{cases}\sum_{l=2}^{N_u}\pk*{\exists_{s' \in E_{u,1}^1, t' \in E_{u,l}^2}\forall_{s \in [s',s'+\frac{S_1}{u^2}], t \in [t',t'+\frac{S_2}{u^2}]}: \begin{array}{ccc}W_1^*(s)>u\\ W_2^*(t)> au\end{array}  } ,& c_2+2c_1 \le 0 \\
\sum_{l=-K_u^{(1)}}^{N_u}\pk*{\exists_{s' \in E_{u,1}^1, t' \in E_{u,l}^2}\forall_{s \in [s',s'+\frac{S_1}{u^2}], t \in [t',t'+\frac{S_2}{u^2}]}:\begin{array}{ccc} W_1^*(s)>u \\ W_2^*(t)> au \end{array} }, & c_2+2c_1 > 0 \end{cases}, $$
which leads us to
$$P_{u,\Delta}^{(1)} \sim \begin{cases}C_{4,\mathcal{P}}^{(1)}(\Delta)C_{4,\mathcal{P}}^{(2)}(\Delta)\frac{\sqrt{3\pi}}{\sqrt{2}} u^{-1}\varphi_{t^*}(u+c_1 ,u+c_2), & c_2+2c_1 \le 0 \\
C_{4,\mathcal{P}}^{(1)}(\Delta)C_{4,\mathcal{P}}^{(2)}(\Delta)\frac{\sqrt{3\pi}}{\sqrt{2}} e^{-2\frac{(\frac{1}{2}c_1+c_2)^2}{3}} \Phi\left(c_2+\frac{1}{2}c_1\right) u^{-1}\varphi_{t^*}(u+c_1 ,u+c_2) , & c_2+2c_1 > 0 \end{cases} $$
as $u \to \IF$ with
$$C_{4,\mathcal{P}}^{(1)}(\Delta)=\int_{\R}\pk*{ \exists_{ s'\in [0,\Delta]} \forall_{s \in [s',s'+S_1]}: W_1(s) - 2s> x} e^{ 2x} dx$$
and
$$C_{4,\mathcal{P}}^{(2)}(\Delta)=\int_{\R} \pk*{ \exists_{ t'\in [0,\Delta]} \forall_{t \in [t',t'+S_2]}: W_2(t) - t> x} e^{ 2x} dx.$$
Similarly, by interchanging the roles of $c_1$ and $c_2$ and using previous results we have that as $u \to \IF$
$$P_{u,\Delta}^{(2)} \sim \begin{cases}C_{4,\mathcal{P}}^{(1)}(\Delta)C_{4,\mathcal{P}}^{(2)}(\Delta)\frac{\sqrt{3\pi}}{\sqrt{2}} u^{-1}\varphi_{t^*}(u+c_1 ,u+c_2), & c_1+2c_2 \le 0 \\
C_{4,\mathcal{P}}^{(1)}(\Delta)C_{4,\mathcal{P}}^{(2)}(\Delta)\frac{\sqrt{3\pi}}{\sqrt{2}} e^{-2\frac{(\frac{1}{2}c_2+c_1)^2}{3}}  \Phi\left(c_1+\frac{1}{2}c_2\right) u^{-1}\varphi_{t^*}(u+c_1 ,u+c_2) , & c_1+2c_2 > 0. \end{cases} $$
Therefore with \nelem{1dFiniteParisian} we can write that
$$\lim_{\Delta \to \IF}\lim_{u \to \IF}\frac{P_{u,\Delta}}{C C_{4,\mathcal{P}}^{(1)}C_{4,\mathcal{P}}^{(2)}u^{-1}\varphi_{t^*}(u+c_1 ,u+c_2)}=1,$$
where
$$C=
\begin{cases}
e^{-2\frac{(\frac{1}{2}c_1+c_2)^2}{3}} \Phi\left(c_2+\frac{1}{2}c_1\right)
+e^{-2\frac{(\frac{1}{2}c_2+c_1)^2}{3}}  \Phi\left(c_1+\frac{1}{2}c_2\right), & c_2>\max(-\frac{1}{2}c_1,-2c_1)\\
e^{-2\frac{(\frac{1}{2}c_1+c_2)^2}{3}} \Phi\left(c_2+\frac{1}{2}c_1\right)
+1, &-\frac{1}{2}c_1<c_2 \le -2c_1 \\
1
+e^{-2\frac{(\frac{1}{2}c_2+c_1)^2}{3}}  \Phi\left(c_1+\frac{1}{2}c_2\right), &-2c_1<c_2 \le -\frac{1}{2}c_1 \\
2 , &c_2\le \min(-\frac{1}{2}c_1,-2c_1).
\end{cases} $$

$D_{u,\Delta}$ is the same as \cite{DHK20}(3.15) or (3.16), exactly in the same way as in the proof of \cite{DHK20}[Thm. 2.2, case (iv)] and hence
\BQNY
\lim_{\Delta \to \IF}\lim_{u \to \IF}\frac{D_{u,\Delta}}{Cu^{-1}\varphi_{t^*}(u+c_1 ,au+c_2)}&=&0.
\EQNY
With that the proof of case (iii) is complete.
\newline
\underline{Case (iv): $\rho<\frac{1}{4a}(1-\sqrt{8a^2+1}).$} From \nelem{optDrift} we have exactly one minimizer of
$q_{\vk a_u(s,t)}^*(s,t)$  on $[0,1]^2$ which is $(s_u,t_u)=(1,\frac{a}{\rho(2a\rho-1)+\frac{c_2-\rho c_1}{u}})$ and for large enough $u$ we have $t_u<1.$
The proof is analogous to case (ii) with
$$F_u:=[1-\frac{\Delta}{u^2},1]\times [t_u-\frac{\log(u)}{u},t_u+\frac{\log(u)}{u}],$$
$$P_{u,\Delta}:=\sum_{l=-N_u}^{N_u}\pk*{\exists_{s' \in E_{u,1}^1, t' \in E_{u,l}^2}\forall_{s \in [s',s'+\frac{S_1}{u^2}]}, t \in [t',t'+\frac{S_2}{u^2}]: W_1^*(s)>u, W_2^*(t)> au  },$$
$$D_{u,\Delta}:=\sum_{l=-N_u}^{N_u}\sum_{m=l+1}^{N_u}\pk*{\exists_{s \in E_{u,1}^1, t_1 \in E_{u,l}^2, t_2 \in E_{u,m}^2}: W_1^*(s)>u, W_2^*(t_1)> au, W_2^*(t_2)> au  }.$$
Using Taylor expansion we have
$$u^2(q_a(k_u,l_u)-q_a(1,t^*))=\tau_1(k-1)\Delta+\tau_4\frac{(l-1)^2\Delta^2}{u^2}+o(\frac{k^2}{u^2})+o(\frac{l^3}{u^4}),$$
where $\tau_1=(1 - 2 a \rho)^2>0, \tau_4=-\frac{ \rho^3 (1 - 2 a \rho)^4}{a (1 - a \rho)}>0.$
Using \nelem{PickandsParisian} and Remark \ref{32} we get
\BQNY
P_{u,\Delta}& \sim & C_{5,\mathcal{P}}^{(1)}(\Delta) C_{5,\mathcal{P}}^{(2)}(\Delta) u^{-2}\varphi_{t_u}(u+c_1 ,au+c_2t^*) \sum_{l=-N_u}^{N_u} e^{-\frac{1}{2}u^2(q_{\vk a}(k_u,l_u)-q_{\vk a}(1,1))}\\
&\sim& C_{5,\mathcal{P}}^{(1)}(\Delta) C_{5,\mathcal{P}}^{(2)}(\Delta) u^{-2}e^{-a\frac{(c_1\rho-c_2)^2}{2\rho(1-a\rho)}} \varphi_{t^*}(u+c_1 ,au+c_2t^*) \sum_{l=-N_u}^{N_u} e^{-\frac{1}{2}u^2(q_{\vk a}(k_u,l_u)-q_{\vk a}(1,1))},\\
\EQNY
where
$$C_{5,\mathcal{P}}^{(1)}(\Delta)=\int_{\R}\pk*{ \exists_{ s'\in [0,\Delta]} \forall_{s \in [s',s'+S_1]}: W_1(s) - \frac{1-a\rho}{1-\rho^2t^*}s> x} e^{ \frac{1-a\rho}{1-\rho^2t^*}x} dx$$
and
$$C_{5,\mathcal{P}}^{(2)}(\Delta)=\int_{\R} \pk*{ \exists_{ t'\in [0,\Delta]} \forall_{t \in [t',t'+S_2]}: W_2(t) - \frac{a}{t^*}t> x} e^{ 2\frac{a}{t^*}x} dx.$$
With \cite{DHK20}[Lemma 3.7] we have
\BQNY
\sum_{l=-N_u}^{N_u} e^{-\frac{1}{2}u^2(q_{\vk a}(k_u,l_u)-q_{\vk a}(1,1))} &=& \frac{u}{\sqrt{\tau_4}\Delta} \sum_{l=1}^{N_u} \frac{\sqrt{\tau_4}\Delta}{u} e^{-\frac{\tau_4}{2}\frac{(l-1)^2\Delta^2}{u^2}}\\
&\sim& \frac{u}{\Delta} \frac{\sqrt{2\pi}}{\tau_4}.
\EQNY
From \nelem{1dFiniteParisian} we have that $\lim_{\Delta \to \IF}C_{5,\mathcal{P}}^{(1)}(\Delta)=C_{5,\mathcal{P}}^{(1)}$ and $\lim_{\Delta \to \IF}\frac{C_{5,\mathcal{P}}^{(2)}(\Delta)}{\Delta}=C_{5,\mathcal{P}}^{(2)}.$ Hence
$$\lim_{\Delta \to \IF}\lim_{u \to \IF}\frac{P_{u,\Delta}}{C_{5,\mathcal{P}}^{(1)}C_{5,\mathcal{P}}^{(2)}\frac{\sqrt{2 \pi}}{\sqrt{\tau_4}} e^{-a\frac{(c_1\rho-c_2)^2}{2\rho(1-a\rho)}} u^{-1}\varphi_{t^*}(u+c_1 ,au+c_2t^*)}=1.$$
$D_{u,\Delta}$ is exactly the same as \cite{DHK20}(3.18) and hence
\BQNY
\lim_{\Delta \to \IF}\lim_{u \to \IF}\frac{D_{u,\Delta}}{Cu^{-1}\varphi_{t^*}(u+c_1 ,au+c_2)}&=&0,
\EQNY
which means that the double sum is negligible. With that the proof of case (iv) is complete.
\newline
\underline{Case (v): $a=1, \rho< A_a.$} According to \nelem{optDrift}, there are two optimal points:
$$(s_u,t_u)= (1,\frac{1}{\rho(2\rho-1)+\frac{c_2-\rho c_1}{u}}), \quad (\bar{s}_u,\bar{t}_u)= (\frac{1}{\rho(2\rho-1)+\frac{c_1-\rho c_2}{u}},1),$$ where for large enough $u$ we have
$$\frac{1}{\rho(2\rho-1)+\frac{c_2-\rho c_1}{u}},\frac{1}{\rho(2\rho-1)+\frac{c_1-\rho c_2}{u}}<1.$$ We can use \eqref{lBP} and \eqref{uBP} with
$$F_u:=[1-\frac{\Delta}{u^2},1]\times[t_u-\frac{\log(u)}{u},t_u+\frac{\log(u)}{u}]\cup [t_u-\frac{\log(u)}{u},t_u+\frac{\log(u)}{u}]\times[1-\frac{\Delta}{u^2},1],$$
\BQNY
P_{u,\Delta}&:=&\pk*{\exists_{(s',t') \in [1-\frac{\Delta}{u^2},1]\times[t_u-\frac{\log(u)}{u},t_u+\frac{\log(u)}{u}]}\forall_{s \in [s',s'+\frac{S_1}{u^2}]}, t \in [t',t'+\frac{S_2}{u^2}]W_1^*(s)>u,W_2^*(t)>u}\\
&&+\pk{\exists_{(s',t') \in [t_u-\frac{\log(u)}{u},t_u+\frac{\log(u)}{u}]\times[1-\frac{\Delta}{u^2},1]}\forall_{s \in [s',s'+\frac{S_1}{u^2}]}, t \in [t',t'+\frac{S_2}{u^2}]W_1^*(s)>u,W_2^*(t)>u}\\
&:=&P_{u,\Delta}^{(1)}+P_{u,\Delta}^{(2)},
\EQNY
\BQNY
D_{u,\Delta}&:=&\pk*{\exists_{(s_1,t_1) \in [1-\frac{\Delta}{u^2},1]\times[t_u-\frac{\log(u)}{u},t_u+\frac{\log(u)}{u}], (s_2,t_2) \in [t_u-\frac{\log(u)}{u},t_u+\frac{\log(u)}{u}]\times[1-\frac{\Delta}{u^2},1]}\begin{array}{ccc}
	W_1^*(s_1)>u \\
    W_1^*(s_2)>u \\
	W_2^*(t_1) > u \\
    W_2^*(t_2) > u
\end{array}}.
\EQNY

Using the same calculations as in case (iv) for $a=1$ we have that
\BQNY
P_{u,\Delta}^{(1)} &\sim& C_{6,\mathcal{P}}^{(1)}C_{6,\mathcal{P}}^{(2)}\frac{\sqrt{2 \pi}}{\sqrt{\tau_4}} u^{-1} \varphi_{t_u}(u+c_1 ,u+c_2t_u)\\
&=& C_{6,\mathcal{P}}^{(1)}C_{6,\mathcal{P}}^{(2)}\frac{\sqrt{2 \pi}}{\sqrt{\tau_4}} e^{-\frac{c_1^2 \rho^2 -2 c_1 c_2 \rho +c_2^2}{2\rho(1-\rho)}}u^{-1} \varphi_{t^*}(u+c_1 ,u+c_2t^*)
\EQNY
where
$$C_{6,\mathcal{P}}^{(1)}(\Delta)=\int_{\R}\pk*{ \exists_{ s'\in [0,\Delta]} \forall_{s \in [s',s'+S_1]}: W_1(s) - \frac{1-\rho}{1-\rho^2t^*}s> x} e^{ \frac{1-\rho}{1-\rho^2t^*}x} dx$$
and
$$C_{6,\mathcal{P}}^{(2)}(\Delta)=\int_{\R} \pk*{ \exists_{ t'\in [0,\Delta]} \forall_{t \in [t',t'+S_2]}: W_2(t) - \frac{1}{t^*}t> x} e^{ \frac{2}{t^*}x} dx.$$
By interchanging the roles of $c_1$ and $c_2$ we can analogously get
\BQNY
P_{u,\Delta}^{(2)} &\sim& C_{6,\mathcal{P}}^{(1)}(\Delta)C_{6,\mathcal{P}}^{(2)}(\Delta)\frac{\sqrt{2 \pi}}{\sqrt{\tau_4}} u^{-1} \frac{1-\rho^2 t^*}{1-\rho} e^{-\frac{c_1^2  -2 c_1 c_2 \rho +c_2^2\rho^2}{2\rho(1-\rho)}}\varphi_{t^*}(u+c_2 ,u+c_1t^*).
\EQNY
From the proof of \cite{DHK20} [Thm 2.2, case (vi)] we have
\begin{enumerate}
\item For $c_1>c_2: \varphi_{t^*}(u+c_1 ,u+c_2t^*)=o(\varphi_{t^*}(u+c_2 ,u+c_1t^*)),$
\item For $c_1<c_2: \varphi_{t^*}(u+c_2 ,u+c_1t^*)=o(\varphi_{t^*}(u+c_1 ,u+c_2t^*)),$
\item For $c_1=c_2: \varphi_{t^*}(u+c_1 ,u+c_2t^*)=\varphi_{t^*}(u+c_2 ,u+c_1t^*).$
\end{enumerate}
With that we obtain that
$$\lim_{\Delta \to \IF}\lim_{u \to \IF}\frac{P_{u,\Delta}}{C_6 C_{6,\mathcal{P}}^{(1)}C_{6,\mathcal{P}}^{(2)}u^{-1}\varphi_{t^*}(u+\min(c_1,c_2) ,u+\max(c_1,c_2)t^*)}=1$$
with
$$ C_6= \begin{cases}
e^{-\frac{( \min(c_1,c_2)\rho- \max(c_1,c_2))^2}{2\rho(1-\rho)}} \frac{\sqrt{2\pi}}{\sqrt{\tau}}, & c_1 \not =c_2 \\	2e^{-\frac{c_2^2(1-\rho)}{2\rho}} \frac{\sqrt{2\pi}}{\sqrt{\tau}}, & c_1=c_2
\end{cases}.$$
$D_{u,\Delta}$ is the same as \cite{DHK20}(3.19) and hence
\BQNY
\lim_{\Delta \to \IF}\lim_{u \to \IF}\frac{D_{u,\Delta}}{Cu^{-1}\varphi_{t^*}(u+c_1 ,au+c_2)}&=&0.
\EQNY
With that the proof of case (v) is complete.
\QED

\section{Appendix}

\prooflem{PickandsParisian}
Let $S_1,S_2>0$. For all the cases we can write:
\BQNY
\lefteqn{\mathcal{P}_{E_{u,k,l},\frac{(S_1,S_2)}{u^2}}(c_1,c_2,u,au)}\\
&=&
\int_{\R^2}\pk*{\exists_{(s',t') \in E} \forall_{s \in [s',s'+S_1], t \in [t',t'+S_2]}:
    \begin{array}{ccc}
	W_1^*(\frac{s}{u^2}+ k_u )>u \\
	W_2^*(\frac{t}{u^2}+ l_u ))> au
	\end{array}
    \Bigg{|}
    \begin{array}{ccc}
	W_1^*(k_u )= u - \frac{x}{u}  \\
	W_2^*(l_u )= au - \frac{y}{u}
	\end{array} } \\	
&\times& u^{-2}\varphi_{k_u,l_u}( u + c_1k_u - \frac{x}{u} ,  au + c_2l_u - \frac{y}{u} ) dxdy\\
&=&
\int_{\R^2}\pk*{\exists_{(s',t') \in E} \forall_{s \in [s',s'+S_1], t \in [t',t'+S_2]}:
    \begin{array}{ccc}
	W_1^*(\frac{s}{u^2}+k_u ) - W_1(k_u ) +W_1(k_u)>u \\
	W_2^*(\frac{t}{u^2}+l_u )- W_2(l_u ) +W_2(l_u)> au
	\end{array}
    \Bigg{|}
    \begin{array}{ccc}
	W_1^*(k_u )= u - \frac{x}{u}  \\
	W_2^*(l_u )= au - \frac{y}{u}
	\end{array}}  \\	
&\times& u^{-2}\varphi_{k_u,l_u}( u + c_1 k_u - \frac{x}{u} ,  au + c_2 l_u - \frac{y}{u}    )dxdy \\
&=&
\int_{\R^2}\pk*{\exists_{(s',t') \in E} \forall_{s \in [s',s'+S_1], t \in [t',t'+S_2]}:
    \begin{array}{ccc}
	\eta_{u,k,l}(s,t)>(x,y)
	\end{array}
    \Bigg{|}
    \begin{array}{ccc}
	W_1^*(k_u )= u - \frac{x}{u}  \\
	W_2^*(l_u )= au - \frac{y}{u}
	\end{array} } \\	
&\times& u^{-2}\varphi_{k_u,l_u}( u + c_1 k_u - \frac{x}{u} ,  au + c_2 l_u - \frac{y}{u}    )dxdy.
\EQNY

From \cite{DHK20} [Lemma 3.3] we have that both for $k_u>l_u$ and $k_u<l_u,$ as $u \to \IF$
\BQN \label{density}
\varphi_{k_u,l_u}( u + c_1 k_u - \frac{x}{u} ,  au + c_2 l_u - \frac{y}{u}    )&\sim& \varphi_{t_u}( u + c_1 ,  au + c_2 t_u) \\
&&\times e^{-\frac{1}{2} u^2(q_{\vk a_u(k_u,l_u)}(k_u,l_u)-q_{\vk a_u(1,t_u)}(1,t_u))}e^{\lambda_1 x + \lambda_2 y}. \nonumber
\EQN
Hence it remains to investigate
$$\int_{\R^2}\pk*{\exists_{(s',t') \in E} \forall_{s \in [s',s'+S_1], t \in [t',t'+S_2]}:
    \begin{array}{ccc}
	\eta_{u,k,l}(s,t)>(x,y)
	\end{array}
    \Bigg{|}
    \begin{array}{ccc}
	W_1^*(k_u )= u - \frac{x}{u}  \\
	W_2^*(l_u )= au - \frac{y}{u}
	\end{array} }  e^{\lambda_1 x + \lambda_2 y} dxdy. $$
Since the interplay between $k_u$ and $l_u$ influences the behaviour of the integrals above, we split the proof into three parts : $k_u=l_u,k_u<l_u,k_u>l_u.$

$(i)$
If $k_u=l_u$, then $\eta^*_{u,k,l,x,y}(s,t)=(\eta^*_{1,u,k,l,x,y}(s),\eta^*_{2,u,k,l,x,y}(t)):=\left(\eta_{u,k,l}(s,t)\Bigg{|}
    \begin{array}{ccc}
	W_1^*(k_u )= u - \frac{x}{u}  \\
	W_2^*(l_u )= au - \frac{y}{u}
	\end{array}\right)$, $s,t \in [-\Delta,0].$ We have that
\[\mathbb{E} \{\eta^*_{u,k,l,x,y}(s,t)\}=-\frac{1}{uk_u}\Bigg{(}\begin{array}{ccc}
	 s(u+c_1k_u-\frac{x}{u})\\
	t(au+c_2k_u-\frac{y}{u})
	\end{array}\Bigg{)}\] and the covariance matrix is equal to
\BQNY
\Sigma_{\left(\eta^*_{u,k,l,x,y}(s,t)\right)}
&=&\Bigg{(}\begin{array}{ccc}
	s & \rho \min(s,t) \\
	\rho \min(s,t) & t
\end{array}\Bigg{)}-u^{-2}\Bigg{(}\begin{array}{ccc}
	s & \rho s \\
	\rho t & t
\end{array}\Bigg{)}
\Bigg{(}\begin{array}{ccc}
	k_u & \rho k_u \\
	\rho k_u & l_u
\end{array}\Bigg{)}^{-1}
\Bigg{(}\begin{array}{ccc}
	s & \rho t \\
	\rho s & t
\end{array}\Bigg{)}\\
&=&\Bigg{(}\begin{array}{ccc}
	s & \rho \min(s,t) \\
	\rho \min(s,t) & t
\end{array}\Bigg{)}-O\left(\frac{\log(u)}{u}\right)\Bigg{(}\begin{array}{ccc}
	s^2 & \rho^2 st \\
	(\rho \min(s,t))^2 & t^2
\end{array}\Bigg{)},s,t\in[0,\Delta].
\EQNY
Similarly
\BQNY
\Sigma_{\left(\eta^*_{u,k,l,x,y}(s_1,t_1)-\eta^*_{u,k,l,x,y}(s_2,t_2)\right)}&=&\left(\begin{array}{ccc}
	|s_1-s_2| & \rho \min(s_1-s_2,t_1-t_2) \\
	\rho \min(s_1-s_2,t_1-t_2) & |t_1-t_2|
	\end{array}\right)\\
&-&O\left(\frac{\log(u)}{u}\right)\left(\begin{array}{ccc}
	|s_1-s_2|^2 & (\rho \min(s_1-s_2,t_1-t_2))^2 \\
	(\rho \min(s_1-s_2,t_1-t_2))^2 & |t_1-t_2|^2
	\end{array}\right).
\EQNY
Together with the continuous mapping theorem we get, as $u \to \infty$
\BQNY
\int_{\R^2}\pk*{\exists_{(s',t') \in E} \forall_{s \in [s',s'+S_1], t \in [t',t'+S_2]}:
    \begin{array}{ccc}
	\eta^*_{u,k,l,x,y}(s,t)>(x,y)
	\end{array}}  e^{\lambda_1 x + \lambda_2 y} dxdy\\
&&
\hspace*{-10cm} \sim \int_{\R^2}\pk*{\exists_{s',t' \in [0,\Delta]} \forall_{s \in [s',s'+S_1], t \in [t',t'+S_2]}:
    \begin{array}{ccc}
	W_1(s)-s>x \\
    W_2(t)-at>y
	\end{array}}  e^{\lambda_1 x + \lambda_2 y} dxdy.
\EQNY
It remains to show that \eqref{Int1} is finite and to justify the use of the dominated convergence theorem, but using \eqref{comparison} we can bound the \eqref{Int1} by replacing the Parisian ruin with simple supremum and the finitness then follows from \cite{DHK20} (3.8).
\newline
$(ii)$
If  $k_u<l_u$, then observe that the increments $W_1(s+k_u u^2) - W_1(k_u u^2),W_2(t+l_u u^2)- W_2(l_u u^2)$ are independent. Hence
\BQNY
\int_{\R^2}\pk*{\exists_{(s',t') \in E} \forall_{s \in [s',s'+S_1], t \in [t',t'+S_2]}:
    \begin{array}{ccc}
	\eta_{u,k,l}(s,t)>(x,y)
	\end{array}
    \Bigg{|}
    \begin{array}{ccc}
	W_1^*(k_u )= u - \frac{x}{u}  \\
	W_2^*(l_u )= au - \frac{y}{u}
	\end{array} }  e^{\lambda_1 x + \lambda_2 y} dxdy \\
&&\hspace*{-15cm} = \int_{\R^2}\pk*{\exists_{s' \in [0,\Delta]}\forall_{s \in [s',s'+S_1]}: \upsilon^*_{1,u,x,y}(s)>x}\\
&&\hspace*{-15cm} \times \pk*{\exists_{t' \in [0,\Delta]}\forall_{t \in [t',t'+S_2]}: \upsilon^*_{2,u,x,y}(t)>y}  e^{\lambda_1 x + \lambda_2 y} dxdy,
\EQNY
where $\upsilon^*_{1,u,x,y}(s):=\left(\eta_{1,u,k}(s)\Bigg{|}
    \begin{array}{ccc}
	W_1^*(k_u )= u - \frac{x}{u}  \\
	W_2^*(l_u )= au - \frac{y}{u}
	\end{array}\right)$ is a Gaussian process with
\[\mathbb{E} \{\upsilon^*_{1,u,x,y}(s)\}=\frac{1}{uk_u}s (u+c_1-\frac{x}{u})- c_1\frac{s}{u},\]
\[Var\left(\upsilon^*_{1,u,x,y}(s)\right)=s-\frac{s^2}{u^2k_u}\] and $\upsilon^*_{2,u,x,y}(t):=\left(\eta_{2,u,l}(t)\Bigg{|}
    \begin{array}{ccc}
	W_1^*(k_u )= u - \frac{x}{u}  \\
	W_2^*(l_u )= au - \frac{y}{u}
	\end{array}\right)$ is a Gaussian process with
\[\mathbb{E}\{\upsilon^*_{2,u,x,y}(t)\}=\frac{1}{u(l_uk_u - \rho^2 k_u^2)}(k_u t (au+c_2-\frac{y}{u})-\rho k_u t (u+c_1-\frac{x}{u}))- c_2\frac{t}{u},\]
\[Var \left(\upsilon^*_{2,u,x,y}(t)\right)=t-\frac{t^2}{u^2(l_u-\rho^2k_u)}.\]
Moreover, for each $0\ge s>t\ge-\Delta$,
$\left(\upsilon^*_{1,u,x,y}(s)-\upsilon^*_{1,u,x,y}(t)\right)$ is Normally distributed with
\[Var \left(\upsilon^*_{1,u,x,y}(s)-\upsilon^*_{1,u,x,y}(t)\right)=(s-t)-\frac{(s-t)^2}{u^2k_u}
\]
while
$\left(\upsilon^*_{2,u,x,y}(s)-\upsilon^*_{2,u,x,y}(t)\right)$ is Normally distributed with
\[Var \left(\upsilon^*_{2,u,x,y}(s)-\upsilon^*_{2,u,x,y}(t)\right)=(s-t)-\frac{(s-t)^2}{u^2(l_u-\rho^2k_u)}.
\]

Hence, using that $Var\left(\upsilon^*_{i,u,x,y}(s)-\upsilon^*_{i,u,x,y}(t)\right)\le 2|s-t|$ for all large enough $u$,\\
$\upsilon^*_{1,u,x,y}(s), s\in[0,\Delta]$ weakly converges in $C[0,\Delta]$ to $W_1(s)-s, s\in[0,\Delta]$ and\\
$\upsilon^*_{2,u,x,y}(t), s\in[0,\Delta]$ weakly converges in $C[0,\Delta]$ to $W_2(s)-\frac{a-\rho}{t^*-\rho^2}s, s\in[0,\Delta].$
It remains to show that \eqref{Int1} is finite and to justify the use of the dominated convergence theorem. As in previous case it falls directly from \eqref{comparison} and \cite{DHK20}(3.8). Combining it with the proven weak convergence and with the dominated convergence theorem, we obtain that
\BQNY
\int_{\R^2}\pk*{\exists_{(s',t') \in E} \forall_{s \in [s',s'+S_1], t \in [t',t'+S_2]}:
    \begin{array}{ccc}
	\eta_{u,k,l}(s,t)>(x,y)
	\end{array}
    \Bigg{|}
    \begin{array}{ccc}
	W_1^*(k_u )= u - \frac{x}{u}  \\
	W_2^*(l_u )= au - \frac{y}{u}
	\end{array} }  e^{\lambda_1 x + \lambda_2 y} dxdy \\
&&\hspace*{-16cm} \sim \int_{\R^2}\pk*{\exists_{s' \in [0,\Delta]}\forall_{s \in [s',s'+S_1]}: W_1(s)-s>x}\\
&&\hspace*{-16cm} \times \pk*{\exists_{t' \in [0,\Delta]}\forall_{t \in [t',t'+S_2]}:  W_2(t)-\frac{a-\rho}{t^*-\rho^2}t>y}  e^{\lambda_1 x + \lambda_2 y} dxdy
\EQNY
\newline
$(iii)$ If $k_u>l_u$, then again the increments $W_1(s+k_u u^2) - W_1(k_u u^2), W_2(t+l_u u^2)- W_2(l_u u^2)$ are independent. Hence
\BQNY
\int_{\R^2}\pk*{\exists_{(s',t') \in E} \forall_{s \in [s',s'+S_1], t \in [t',t'+S_2]}:
    \begin{array}{ccc}
	\eta_{u,k,l}(s,t)>(x,y)
	\end{array}
    \Bigg{|}
    \begin{array}{ccc}
	W_1^*(k_u )= u - \frac{x}{u}  \\
	W_2^*(l_u )= au - \frac{y}{u}
	\end{array} }  e^{\lambda_1 x + \lambda_2 y} dxdy \\
&&\hspace*{-15cm} =\int_{\R^2}\pk*{\exists_{s' \in [0,\Delta]}\forall_{s \in [s',s'+S_1]}: \nu^*_{1,u,x,y}(s)>x}\\
&&\hspace*{-15cm} \times \pk*{\exists_{t' \in [0,\Delta]}\forall_{t \in [t',t'+S_2]}: \nu^*_{2,u,x,y}(t)>y}  e^{\lambda_1 x + \lambda_2 y} dxdy,
\EQNY

where $\nu^*_{1,u,x,y}(s):=\left(\eta_{1,u}(s)\Bigg{|}
    \begin{array}{ccc}
	W_1^*(k_u )= u - \frac{x}{u}  \\
	W_2^*(l_u )= au - \frac{y}{u}
	\end{array}\right)$ is a Gaussian process with
\[\mathbb{E} \{\nu^*_{1,u,x,y}(s)\}=\frac{1}{u(l_uk_u - \rho^2 l_u^2)}(s l_u (u+c_1-\frac{x}{u})-\rho s l_u (au+c_2-\frac{y}{u}))- c_1\frac{s}{u},\]
\[Var\left(\nu^*_{1,u,x,y}(s)\right)=s-\frac{s^2}{u^2(k_u-\rho^2l_u)}\] and $\nu^*_{2,u,x,y}(t):=\left(\eta_{2,u}(t)\Bigg{|}
    \begin{array}{ccc}
	W_1^*(k_u )= u - \frac{x}{u}  \\
	W_2^*(l_u )= au - \frac{y}{u}
	\end{array}\right)$ is a Gaussian process with
\[\mathbb{E} \{\nu^*_{2,u,x,y}(t)\}=\frac{1}{ul_u}t(au+c_2-\frac{y}{u})- c_2\frac{t}{u},\]
\[Var \left(\nu^*_{2,u,x,y}(t)\right)=t-\frac{t^2}{u^2l_u}.\]
Moreover, for each $0\ge s>t\ge-\Delta$,
$\left(\nu^*_{1,u,x,y}(s)-\nu^*_{1,u,x,y}(t)\right)$ is normally distributed with
\[Var \left(\nu^*_{1,u,x,y}(s)-\nu^*_{1,u,x,y}(t)\right)=(s-t)-\frac{(s-t)^2}{u^2(k_u-\rho^2l_u)}
\]
and
$\left(\nu^*_{2,u,x,y}(s)-\nu^*_{2,u,x,y}(t)\right)$ is normally distributed with
\[Var \left(\nu^*_{2,u,x,y}(s)-\nu^*_{2,u,x,y}(t)\right)=(s-t)-\frac{(s-t)^2}{u^2l_u}.
\]

Hence, using that $Var \left(\nu^*_{i,u,x,y}(s)-\nu^*_{i,u,x,y}(t)\right)\le 2|s-t|$ for all u large enough,\\
$\nu^*_{1,u,x,y}(s), s\in[0,\Delta]$ weakly converges in $C[0,\Delta]$ to $W_1(s)-\frac{1-a\rho}{1-\rho^2t^*}s, s\in[0,\Delta]$ and\\
$\nu^*_{2,u,x,y}(t), t\in[0,\Delta]$ weakly converges in $C[0,\Delta]$ to $W_2(t)-\frac{a}{t^*}t, t\in[0,\Delta]$.
This leads to
\BQNY
\int_{\R^2}\pk*{\exists_{(s',t') \in E} \forall_{s \in [s',s'+S_1], t \in [t',t'+S_2]}:
    \begin{array}{ccc}
	\eta_{u,k,l}(s,t)>(x,y)
	\end{array}
    \Bigg{|}
    \begin{array}{ccc}
	W_1^*(k_u )= u - \frac{x}{u}  \\
	W_2^*(l_u )= au - \frac{y}{u}
	\end{array} }  e^{\lambda_1 x + \lambda_2 y} dxdy \\
&&\hspace*{-14cm} \sim \int_{\R^2}\pk*{\exists_{s' \in [0,\Delta]}\forall_{s \in [s',s'+S_1]}: W_1(s)-\frac{1-a\rho}{1-\rho^2t^*}s>x} \\
&&\hspace*{-14cm} \times \pk*{\exists_{t' \in [0,\Delta]}\forall_{t \in [t',t'+S_2]}:  W_2(t)-\frac{a}{t^*}t>y}  e^{\lambda_1 x + \lambda_2 y} dxdy.
\EQNY
The finitness of \eqref{Int1} and the application of the dominated convergence theorem can be proven identically as in the previous case. This completes the proof.
\QED

\bibliographystyle{ieeetr}

\bibliography{queue2d}
\end{document}